\documentclass{amsart}[12pt]

\newtheorem{lem}{Lemma}
\newtheorem{thm}{Theorem}
\newtheorem{cor}{Corollary}

\newtheorem{pro}{Proposition}

\newcommand{\di}{\displaystyle}

\usepackage[ marginparsep=5mm,heightrounded,centering,margin=3cm]{geometry}

\newenvironment{defn}[1][Definition]{\begin{trivlist}
\item[\hskip \labelsep {\bfseries #1}]}{\end{trivlist}}

\usepackage{enumerate}
\usepackage{amssymb,amsmath}

\usepackage[numbers]{natbib}

\title{On growth of metabelian Lie algebras}
\author{D\.{I}lber Ko\c{c}ak}
\begin {document}

\begin{abstract} 

For any integer $d\geq 1$  we construct examples of finitely presented 
 algebras with intermediate growth of type $[e^{n^{d/(d+1)}}]$. We produce these examples 
 by computing the growth types of some finitely presented metabelian Lie algebras.
 \end{abstract}
\maketitle

\begin{section}{Introduction}
Let $A$ be an (not necessarily associative) algebra  over a field $k$ generated by a finite set $X$, 
 and $X^n$ 
denote the subspace of $A$ spanned by all monomials on $X$ of length at most
$n$. The 
\textit{growth function} of $A$ with respect to $X$ is defined by
$$\gamma_{X,A}(n)=dim_k(X^n)$$
This function depends on the choice of the generating set $X$. To remove this dependence 
the following 
equivalence relation is introduced:
Let $f(n)$ and $g(n)$ be two monotone functions on  $ \mathbb{N}$. We write $ f\precsim g$ if there exists 
 a constant $C\in \mathbb{N}$ such that $f(n)\leq Cg(Cn)$ for all natural numbers $n$. If $ f\precsim g$
 and $ g\precsim f$, we say $f$ and $g$ are equivalent and denote this by $f\sim g$. The
 equivalence class containing $f$ is denoted  by $[f]$ and is called the 
 \textit{growth rate of f}.
Set $[f]\leq [g]$ if and only if $f\precsim g$.
\par The growth rate $[2^n]$ is called \textit{exponential} and a growth rate strictly less than
exponential
is called \textit{subexponential}. A subexponential growth which is greater than $[n^d]$
for any $d$ is called \textit{intermediate}.
\par The growth rate is a widely studied invariant for
finitely generated algebraic structures such as groups, semigroups and algebras.
 The notion of growth for groups was introduced by Schwarz 
\cite{schwarz55} and independently by Milnor \cite{milnor68}. 
The study of growth of algebras dates back to the papers by Gelfand and Kirillov, 
\cite{GK661,GK66}.
A theorem of Milnor and Wolf \cite{milnor68, Wolf68} states that any solvable group
has a polynomial growth if it is virtually nilpotent, otherwise 
it has exponential growth. The description of groups of 
polynomial growth was obtained by Gromov in his celebrated work \cite{gromov81}. 
He proved that every
finitely generated group of polynomial growth is virtually nilpotent. 
The situation for algebras is
different from the case of groups. M. Smith \cite{smith76} showed that there exists an 
infinite dimensional solvable Lie algebra $L$ whose universal enveloping algebra $U(L)$ has
intermediate growth. Later, in \cite{lichtman84}, Lichtman proved that the universal envelope
of an arbitrary finitely generated infinite dimensional virtually solvable Lie algebra has 
intermediate growth. In \cite{LiUf95}, Lichtman and Ufnarovski showed that the growth 
rate of
a finitely generated free solvable Lie algebra of derived length $k>3$ and its universal
envelope are almost exponential (this means that it is less than exponential growth $[2^n]$
but greater than the growth $[2^{n^{\alpha}}]$ for any $\alpha< 1$).
\par The first examples of finitely generated groups of intermediate growth were constructed by 
Grigorchuk \cite{Gri83, Gri84}. It is still an open problem whether there exists  finitely 
presented groups of intermediate growth. 
In contrast, there are examples of
finitely presented algebras of intermediate growth.
The first such example is the universal enveloping algebra of the
Lie algebra $W$ with basis 
$\{w_{-1},w_0,w_1,w_2,\dots\}$ and brackets defined by $[w_i,w_j]=(i-j)w_{i+j}$. $W$ is a 
subalgebra 
of the 
generalized Witt algebra $W_{\mathbb{Z}}$ (see \cite[p.206]{as74} for definitions).
It was proven in \cite{stewart75} that $W$ has a finite presentation with two generators and 
six relations.  It is also a graded algebra with generators of degree $-1$ and $2$. 
Since $W$ has linear
growth, its universal enveloping algebra has growth $[e^{\sqrt{n}}]$ and it is finitely presented. Known examples of  finitely presented algebras of intermediate growth have growth of 
type $[e^{\sqrt{n}}]$.
In this note we present  examples of finitely presented associative algebras of
intermediate growth having different growth types. Specifically, our main result is the following:
\begin{thm}
 For any positive integer $d$, there exists a finitely presented associative algebra
 with intermediate  growth of type $[e^{n^{d/d+1}}]$.
\end{thm}

The steps in proving the theorem are as follows:
In \cite{bau77}, Baumslag established the fact that every finitely generated
metabelian Lie algebra can be embedded in a finitely presented metabelian Lie algebra.
Using ideas of \cite{bau77} (and clarifying some arguments thereof), we embed the free metabelian Lie algebra $M$ (with $d$ generators) into a finitely 
presented metabelian Lie algebra $W^+$.   Next we show that $W^+$ has polynomial growth of type $[n^d]$.
Finally, considering the universal  enveloping algebra of $W^+$ we obtain a finitely  presented associative algebra of growth type $[e^{n^{d/d+1}}]$.

\section{Growth of a finitely generated free metabelian Lie algebra}
Let $k$ be a field and $L$  a Lie algebra over $k$ generated by a finite set $X$. Elements of $X$ are monomials of length $1$. Inductively,  a monomial of length $n$ is an element of the form 
$[u,v]$, where $u$ is a monomial of length $i<n$ and $v$ is monomial of length $n-i$. Every element of $L$
is a linear combination of monomials.    
  If 
$a_1,\dots,a_n \in X$ then $[a_1,\dots ,a_n]$ is defined inductively by 
$$[a_1,\dots,a_n]=[[a_1,\dots,a_{n-1}],a_n]\;\text{for}\; n>2$$
Monomials of the form $[a_1,\dots,a_n]$ are called \textit{left-normed}.
We will frequently use the following simple lemmas in the remainder of this note.

\begin{lem}\label{l0}
 Let $x,\;y,\;z$ be elements of a Lie algebra and $[x,y]=0$. Then the following relations hold:
 $$[x,[y,z]]=[y,[x,z]]$$
 $$[x,z,y]=[y,z,x].$$
\end{lem}
\begin{proof}

Direct consequence of Jacobi identity.
\end{proof}

\begin{lem}\label{l00}
Any element of a Lie algebra can be written as a linear combination of left-normed monomials.
\end{lem}
\begin{proof}
By induction on the length of monomials.
\end{proof}

\begin{defn}
 A Lie algebra $L$ is called  \textit{solvable of derived length $n$} if 
 $$L^{(n)}=0\;\text{and}\;L^{(n-1)}\neq 0$$
 where $L^{(m+1)}=[L^{(m)},L^{(m)}]$
 and $L^{(0)}=L$. We also denote $L^{(1)}=[L,L]$ by $L^{\prime}$ and call it the
 \textit{commutator} of $L$. A solvable Lie algebra of derived length $2$ is called 
 \textit{metabelian}.
\end{defn}
Let $X=\{x_1,\dots,x_d\}$  be a finite set and $L_X$ be the free Lie algebra generated by $X$. 
$M=L_X/L_{X}^{(2)}$ is the free metabelian Lie algebra generated by $X$. The following 
proposition can be found in \cite{bokut63}.
\begin{pro}\label{p1}
 Let $M$ be a free metabelian Lie algebra over a field $k$
 with the generating set $X=\{x_1,\dots,x_d\}$
 and $x_1<\dots <x_d$ be an order on $X$. If $\mathcal{B}$ is the set of left-normed monomials of the form
 $$[a_0,a_1,\dots,a_{n-1}]$$
 where $a_0>a_1\leq a_2\leq\dots \leq a_{n-1}$ for $a_i\in X$ and $n\geq 1$.
 Then $\mathcal{B}$ forms a basis for $M$.
\end{pro}
\begin{proof} 
 Let $M_n$ denote the subspace of $M$ spanned by all left-normed monomials of length $n$ in $M$
 ($M_0=k$ and $M_1$ is the subspace spanned by $X$). Since the relations $[x,y]=-[y,x]$ for 
 $x,y\in M$ and Jacobi identity are homogeneous,
 $ M_k\cap M_l= \emptyset $ for $k\neq l$. By Lemma \ref{l00}, we get $M= \bigoplus_{n=0}^{\infty} M_n$
  and $\mathcal{B}=\bigsqcup_{n=1}^{\infty}\mathcal{B}_n$ where
 $\mathcal{B}_n$ denotes the set of monomials of length $n$ in $\mathcal{B}$. Hence, it is enough to check that $\mathcal{B}_n$ is a basis of $M_n$ for any $n\geq 1$
i.e.,
\begin{itemize}
 \item[(i)] Any element of $M_n$ can be written as a linear combination of the elements of $B_n$.
 \item [(ii)] Elements of $\mathcal{B}_n$ are linearly independent.
 \end{itemize}
 For $n=1$, $\mathcal{B}_1=X$, so $(i)$ and $(ii)$ hold. 
\\Assume $n=2$. By the anti-symmetry of the Lie bracket
 $$[x,y]=-[y,x]$$
 for any $x,y\in X$.
 So $\mathcal{B}_2$ spans $M_2$.
 \\Assume $n=3$ and $y_1,y_2,y_3\in X$ such that $y_1\leq y_2\leq y_3$. There are at most 6 
 different monomials of length 3 containing $y_1,y_2,y_3$. They are
 $$m_1=[y_1,y_2,y_3]$$
 $$m_2=[y_1,y_3,y_2]$$
 $$m_3=[y_2,y_1,y_3]$$
 $$m_4=[y_2,y_3,y_1]$$
 $$m_5=[y_3,y_1,y_2]$$
 $$m_6=[y_3,y_2,y_1]$$
 $m_3$ and $m_5$ are either $0$ or the elements of $\mathcal{B}_3$ (For example, if $y_1=y_3$
 , $m_5=0$). 
 $$m_1=[y_1,y_2,y_3]=-[y_2,y_1,y_3]=-m_3$$
 $$m_2=[y_1,y_3,y_2]=-[y_3,y_1,y_2]=-m_5$$
and by Jacobi identity,
 $$ \di \begin{array}{lll}
       \di  m_4&=&[y_2,y_3,y_1]\\
                                  &=&-[y_1,y_2,y_3]-[y_3,y_1,y_2]\\
                                  
                                  &=&[y_2,y_1,y_3]-[y_3,y_1,y_2]\\
                                  &=&m_3-m_5
  \end{array}             
$$
and,
 $$ m_6=-m_4=m_5-m_3  $$
 Hence, $\mathcal{B}_3=\{[x_{i_0},x_{i_1},x_{i_2}]\mid x_{i_0}>x_{i_1}\leq x_{i_2},\;x_{i_j}\in X\}$ spans $M_3$.
 \\Now assume that $\mathcal{B}_k$ spans $M_k$ for any $k\in\{2,3,\dots,n\}$. Let 
 $a=[a_0,a_1,\dots,a_n]$ be an element of length $n+1$ in $M$ where $a_i\in X$.
 By the assumption $[a_0,a_1,\dots,a_{n-1}]$ can be written as a linear combination of the
 elements of $\mathcal{B}_n$. So $a$ is a linear combination of the elements of
 the form $[a_{i_0},a_{i_1},\dots,a_{i_{n-1}},a_n]$ where $i_0,i_1,\dots,i_{n-1}$ is
 a permutation
 of $0,1,\dots,n-1$ and $a_{i_0}>a_{i_1}\leq \dots \leq a_{i_{n-1}}$. 
 If $a_n \geq a_i$ for any $i\in \{0,\dots, n-1\}$, then 
 $[a_{i_0},a_{i_1},\dots,a_{i_{n-1}},a_n]$ is in $\mathcal{B}_{n+1}$.
 Otherwise there exists $j\in\{0,\dots,n-1\}$ such that
 $a_{i_j}$
 is the smallest element satisfying $a_{i_j}>a_n$. If $j=0$ then again 
  $[a_{i_0},a_{i_1},\dots,a_{i_{n-1}},a_n]$ is in $\mathcal{B}_{n+1}$. If $j>0$, we apply
  Jacobi identity to  $[a_{i_0},a_{i_1},\dots,a_{i_{n-1}},a_n]$,
  
  $$ \di \begin{array}{lll}
       \di [a_{i_0},a_{i_1},\dots,a_{i_{n-1}},a_n] &=&
       -([[a_{n},[a_{i_0},\dots,a_{i_{n-2}}]],a_{n-1}]+[[a_{n-1},a_n],[a_{i_0},
       \dots,a_{i_{n-2}}]])\\
                                  &=&[a_{i_0},a_{i_1},\dots,a_{i_{n-2}},a_n,a_{n-1}]

  \end{array}             
$$
For $j\geq 1$, we can repeat this $n-j-1$ times and get
\begin{equation}\label{a}
       \di [a_{i_0},a_{i_1},\dots,a_{i_{n-1}},a_n] =
       [a_{i_0},\dots,a_{i_j},a_n,a_{i_{j+1}},\dots,a_{n-1}]
       \end{equation}

If $j\geq 2$, we apply the identity one more time and get
$$ \di \begin{array}{lll}
       \di [a_{i_0},a_{i_1},\dots,a_{i_{n-1}},a_n] &=&
       [a_{i_0},\dots,a_{i_{j-1}},a_n,a_{i_j},a_{i_{j+1}},\dots,a_{n-1}]
                                  
  \end{array}             
$$
and it is a monomial in $\mathcal{B}_{n+1}$. If $j=1$,
$$[a_{i_0},a_{i_1},a_n,a_{i_2},\dots,a_{i_{n-1}}]$$
and since $[a_{i_0},a_{i_1},a_n]$ is an element of the vector space spanned by $\mathcal{B}_3$,
$a$ is a linear combination of monomials in $\mathcal{B}_{n+1}$.
To complete the proof , we also need to verify $(ii)$ for any $n\geq 2$:
\\Since $M_i\cap M_j$ for any $i,j$, $\di \sum_{i=1}^{m}c_ip_i(x_1,\dots,x_d)=0$ where $c_i\in k$ and
$p_i(x_1,\dots,x_d)\in M_i$ implies that either $k_i=0$ or $p_i(x_1,\dots,x_d)=0$. If 
$ p_n(x_1,\dots,x_d)=\sum_{j=1}^{l}c_j[x_{j_1},\dots,x_{j_n}]  $ for $c_j\in k $, $x_{j_t}\in X$, $j\in \{1,\dots,l\}$,
$t\in \{1,\dots,n\}$ then the indices $j_1,\dots,j_n$ are permutations of $1_1,\dots,1_n$ for any 
$j\in \{1,\dots,l\}$.
For $n=2$, $[x,y]\in \mathcal{B}_2$ implies $[y,x]\notin \mathcal{B}_2$. So $\mathcal{B}_2$ forms a basis
 for $M_2$.
 \\For $n=3$, let $m_1,\dots,m_6$
be monomials of length $3$ as we defined before.  Among them only $m_3$ and $m_5$ are in $\mathcal{B}_3$, and
$m_3\neq cm_5$ for any $c\in k$. Hence, $(ii)$ holds for $n=3$. 
\\Now assume that for $3\leq i\leq n$, the elements of $\mathcal{B}_i$ are linearly independent and let 
$[a_0,a_1,\dots,a_n]\in \mathcal{B}_{n+1}$. Assume that for $i\in \{1,\dots, m\}$, $m\in\mathbb{N}$,
there are indices $i_0,\dots,i_n$ which are 
permutations of $0,\dots,n$ such that 
$$[a_{i_0},\dots,a_{i_n}]\in \mathcal{B}_{n+1}\;\text{for any}\; i\in \{1,\dots,m\}$$
and 
$$[a_0,a_1,\dots,a_n]=\sum_{i=1}^{m}c_i[a_{i_0},\dots,a_{i_n}]\;\text{for some}\;c_i\in k$$

Choose the smallest element $a_j$ among $\{a_0,\dots,a_n\}$. It is clear that $j\neq 0$ and if $j=1$ then there
exists $k>1$ such that $a_1=a_k$. So we can assume that $j\geq 2$. Similarly, for any $i\in \{1,\dots,m\}$,
there exists $k_i\geq2$ such that $a_j=a_{i_{k_i}}$. By using equation \ref{a}, we can replace $a_j$ and $a_n$
in $[a_0,\dots,a_n]$ and $a_{i_{k_i}}$ and $a_{i_n}$ in all monomials $[a_{i_0},\dots,a_{i_n}]$, 
$i\in \{1,\dots,m\}$. We get

$$[b,a_j]=\sum_{i=1}^{m}c_i[b_i,a_j]$$
for $b,b_1,\dots,b_m\in M_n$. Any element of $M_n$ can be written as a linear combination of elements of
$\mathcal{B}_n$ but this contradicts that $\mathcal{B}_n$ is linearly independent.
Hence, we conclude that $\mathcal{B}$ is a basis of $M$.

\end{proof}
\begin{cor}\label{c1}
 Let $M$ be a free metabelian Lie algebra over a field $k$ with generating set 
 $X=\{x_1,\dots,x_d\}$. Then $M$
 has polynomial growth of degree $d$.
 
\end{cor}
\begin{proof}
 Let $x_1<\dots <x_d$ be the order on $X$.
The growth function of $M$
 is
 $$\di \gamma_{X,M}(n)=dim_k(X^n)=\sum_{k=1}^{n}|\mathcal{B}_k                                                |$$
 where $\mathcal{B}_k$ denotes the set of basis elements of length $k$ as in Proposition 1.
 For $m>2$, consider $\mathcal{B}_m$. The elements of $\mathcal{B}_m$ are of the form 
 $$a=[a_0,a_1,\dots,a_{m-1}]$$
 where $a_i\in X$ and $a_0>a_1\leq a_2\leq \dots \leq a_{m-1}$.
 \\If $a_1=x_j$ for some $j\in {1,\dots ,d-1}$ then for $i\in {1,\dots,m-1}$,
 $a_i\in \{x_j,\dots,x_d\}$ and $a_0\in \{x_{j+1},\dots,x_d\}$. So for fixed $a_1=x_j$, the 
 number of basis elements of length $m$ is
 $$(d-j)\binom{m-2+d-j}{d-j}$$
 Hence,
 $$
  \begin{array}{lll}
  |B_m|&=&\di\sum_{j=1}^{d-1}(d-j)\binom{m-2+d-j}{d-j}\\
 & &\\
 &=& \di \sum_{j=1}^{d-1}\frac{(m-2+d-1)+\dots + (m-2+1)}{(d-j-1)!}\\
 & &\\
 &\sim&\di \sum_{i=0}^{d-2}\frac{1}{i!}(m-1)^{d-1}\\
 & &\\
 &\di \sim&m^{d-1}
 \end{array}$$
  Since
 $$\di \gamma_{X,M}(m)=dim_k(X^m)=d+\binom{d}{2}+\sum_{k=3}^{m}|\mathcal{B}_k|,$$ 
 $M$ has polynomial growth of degree $d$.
\end{proof}
\section{Wreath Product of two abelian Lie algebras}
Let $T$ be a Lie algebra over a field $k$ of characteristic $p\neq 2$ and $B$ a $k$-module. $B$
is called a \textit{right $T$-module}, if there exists a $k$-bilinear map 
$B\times T\rightarrow B$ satisfying the following :
$$b[t_1,t_2]=(bt_1)t_2-(bt_2)t_1$$
for $b\in B$, $t_1,t_2\in T$.
\\For a Lie algebra $T$ and a right $T$- module $B$, we define the split extension $W=B]T$ of $B$ by
$T$ as follows: As a vector space $W=B\oplus T$ is 
the direct sum of $B$ and $T$, and  the Lie operation on $W$ is defined as
$$[b_1+t_1,b_2+t_2]= (b_1t_2-b_2t_1)+[t_1,t_2]$$
for $b_1,b_2\in B$ and $t_1,t_2\in T$, so $W$ is a Lie algebra over $k$.
\\Here, we consider a special case of this construction that can be found in 
\cite{bau77, lichtman84,bah87}. Suppose $A$ and $T$ are finite dimensional abelian Lie algebras
over a field $k$ of characteristic $p\neq 2$. 
Let 
$\{a_1,\dots,a_m\}$ and $\{t_1,\dots,t_n\}$ be bases of $A$ and $T$, respectively. 
Let $U=U(T)$ be the universal enveloping algebra of $T$
and $B$ be a free right $U$-module with the module basis $\{a_1,\dots,a_m\}$. $B$ can be also viewed
as a Lie algebra module for $T$ where the
action of $T$ on $B$ is the action of a subset of $U$ on the $U$-module $B$
and hence we can form the split extension $W=B]T$. $A$ and $T$ are
Lie subalgebras of $W$, $B$ is the ideal generated of $W$ generated by $\{a_1,\dots,a_m\}$ and it is called 
the \textit{base ideal} of $W$. $W$ is a Lie algebra generated
by $A$ and $T$. It is termed the \textit{wreath product} of the Lie algebras of $A$ and $T$ and denoted by
$$W=A \wr T$$
(For the general definition of the wreath product of Lie algebras  see \cite{shmelkin73,bah87}).
As a vector space $W=B\oplus T$ is the direct sum of its abelian ideal $B$ and abelian Lie subalgebra $T$.
\begin{lem}\label{l1}
  \cite{bau77} Suppose that the Lie algebra $W$ over $k$ is the direct sum of $B$ and $T$ where 
    $B$ is an abelian ideal generated by $\{a_1,\dots,a_m\}$ and $T$ is an abelian Lie algebra
    generated by $\{t_1,\dots,t_n\}$:
    
 $$W=B\oplus T.$$
 Then 
 $B$ is a Lie algebra spanned by 
 $$\{[a_l,t_{j_1},\dots,t_{j_s}]
 \mid \; l \in \{1,\dots,m\}\;\text{and}\; j_1,\dots,j_2\in \{1,2,\dots n\}\}.$$
 \end{lem}
\begin{proof}
 Firstly, we show that $W=B\oplus T$ is metabelian:
 Let $w,w'$ be elements of $W$,
 $$w=b+t\;\text{and}\;w'=b'+t',\;\; b,b'\in B\; t,t'\in T. $$
 Then,
 $$\di \begin{array}{lll}
 
 [w,w']&=&[b+t,b'+t']\\
 &=&[b,b'+t']+[t,b'+t']\\
&=&[b,b']+[b,t']+[t,b']+[t,t'] \\
&=&[b,t']-[b',t] \in B 

 \end{array} $$
 this implies that $W'\subset B$. Since B is abelian, $W^{(2)}=0$. In W,
all the elements 
 can be written as linear combinations of left-normed commutators with elements from
 $\{a_1,\dots,a_m,t_1,\dots,t_m\}$.
 \\Consider a commutator  $x=[x_1,\dots,x_s]$ for $s\geq 2$ and 
 $x_i \in \{a_1,\dots,a_m,t_1,\dots,t_m\}$. If $x\neq 0$ then there is exactly one
 $j$ such that $x_j\in \{a_1,\dots,a_m\}$ (In particular $j=1$ or $j=2$, otherwise $[x_1,x_2]=0$):
If $x_i\in\{t_1,\dots,t_n\}$ for all $i\in \{1,\dots,s\}$, then $x=0$.
If there exist $x_k,x_l \in \{a_1,\dots,a_m\}$ for some $k\neq l$
then $x=[x_1,\dots,x_k,\dots,x_{l-1},x_l,\dots x_s]=0$ since $[x_1,\dots,x_k,\dots,x_{l-1}],x_l \in B$
and their product is equal to $0$.
So $W$ has a basis which is a subset of the following set:
$$\{t_1,\dots,t_d\}\cup \{[a_l,t_{j_1},\dots,t_{j_s}])
 \mid \; l \in \{1,\dots,m\}\;\text{and}\; j_1,\dots,j_2\in \{1,2,\dots n\}\}.$$
 Since $B\cap T=\{0\} $ ,
 $$B\leq span([a_l,t_{j_1},\dots,t_{j_s}]
 \mid \; l \in \{1,\dots,m\}\;\text{and}\; j_1,\dots,j_2\in \{1,2,\dots n\}) $$ and we have 
 $[a_l,t_{j_1},\dots,t_{j_s}]\in B$ for l 
 $\in \{1,\dots,m\}\;\text{and}\; j_1,\dots,j_2\in \{1,2,\dots n\}$. Hence 
  $$B=span([a_l,t_{j_1},\dots,t_{j_s}]
 \mid \; l \in \{1,\dots,m\}\;\text{and}\; j_1,\dots,j_2\in \{1,2,\dots n\}).$$
 
\end{proof}

\begin{cor}\label{c0}
 $W$ can be presented by the generators
$$a_1,\dots,a_m,t_1,\dots,t_n$$
and the following relations:
$$[t_i,t_j]=0$$
for $1\leq i,j\leq n$,
 $$[[a_k,t_{i_1},\dots,t_{i_r}],[a_l,t_{j_1},\dots,t_{j_s}]]=0$$
  for $1\leq k,l\leq m,\;\{i_1,\dots,i_r,j_1,\dots,j_s\}\subset \{1,2,\dots,n\},\;r\geq 0
  ,\;s\geq 0$.
\end{cor}

 The next lemma follows from a theorem of Lewin \cite{lewin74} and it is reformulated in \cite{bau77} as:
 \begin{lem}\label{l3}
  Let $F$ be a finitely generated free Lie algebra 
  and $R$ an ideal of $F$. Then $F/R^{\prime}$ can be embedded in $W=A\wr( F/R)$
  where $A$ is a finite dimensional abelian Lie algebra.
 \end{lem}

\begin{cor}\label{c2}
 Let $M$ be a finitely generated free metabelian Lie algebra. Then there exist finite-dimensional abelian 
 Lie algebras $A$ and $T$ such that $M$ can be embedded in $W=A \wr T$.
\end{cor}

\begin{proof}
Let $R=L_X^{\prime}$ be the commutator of the free Lie algebra $L_X$ generated by $X$. Then by Lemma \ref{l3}, $M=L_X/L_X^{(2)}$ can be embedded 
in $ \di W=A\wr T$ where $T= L_X/L_X^{\prime}$ and $A$ is a finite dimensional abelian Lie algebra.
\end{proof}
\section{Finitely presented metabelian Lie algebras }
Let $A$ and $T$ be finite dimensional abelian Lie algebras over a field $k$ of characteristic $p\neq 2$ and
$W$ the wreath product of $A$ and $T$ as we defined in the previous section. In \cite{bau77}, Baumslag showed that
$W$ can be embedded in a finitely presented metabelian Lie algebra $W^+$. The construction of $W^+$ is as follows:
\\Let $\{a_1,\dots a_m\}$ and $\{t_1,\dots,t_n\}$ be bases of $A$ and $T$, respectively.
 Then the universal enveloping algebra $U$ of $T$ is the associative 
 k-algebra $k[t_1,\dots, t_n]$ of polynomials with variables $t_1,\dots,t_n$ over $k$. Furthermore, let $B$ be the free right $U$-module
 with module basis $\{a_1,\dots,a_m\}$. It is well known that $U$ can be turned into a Lie algebra
 by defining a new multiplication in $U$ by
 $$[u,v]=uv-vu$$
 $U$ is simply an infinite dimensional abelian Lie algebra and $T$ is a finite dimensional subalgebra of
 $U$. To get a finitely presented metabelian Lie algebra $W^+$, we consider a subalgebra $T^+$ of the Lie
 algebra $U$ properly containing $T$ as a subalgebra. Let $T^+$ be the subalgebra generated by 
 $\{t_1,\dots,t_n,u_1,\dots,u_n\}$ where 
 $$u_i=t_i^2\;\;\text{for}\;i\in\{1,\dots,n\}$$
 $T^+$ is a $2n$-dimensional abelian Lie algebra and we define $W^+$ as the wreath product of $A$ and $T^+$.
 $$W^+=A\wr T^+$$
 
 \begin{lem}\cite[Lemma 6]{bau77}\label{l2} $W^+$ can be presented on the generators
 $$a_1,\dots,a_m,t_1,\dots,t_n,u_1,\dots,u_n,$$
 subject to the relations
 $$[[a_k,t_{i_1},\dots,t_{i_r}],[a_l,t_{j_1},\dots,t_{j_s}]]=0,$$
$(1\leq k\leq m,\;1\leq l\leq m, i_1,\dots,i_r,j_1,\dots,j_r \in \{1,2,\dots,n\})$,
 $$[t_i,t_j]=[t_i,u_j]=[u_i,u_j]=0\;\;\;(1\leq i\leq n,\; 1\leq j\leq n),$$
 $$[a_k,u_l]=[a_k,t_l,t_l]\;\;\;(1\leq k\leq n,\; 1\leq l\leq n).$$
  
 \end{lem}
 \begin{proof}By the construction of $W^+$, we have the following relations
 \begin{equation}\label{eq:1}
   [a_k,u_i]=a_kt_i^2=[a_k,t_i,t_i],
 \end{equation}
 and, since $T^+$ is abelian
 \begin{equation}\label{eq:2}
  [t_i,t_j]=[t_i,u_j]=[u_i,u_j]=0
 \end{equation}

  for any $1\leq k\leq m,\;1\leq i,j\leq n.$ It follows from Lemma \ref{l0} and \eqref{eq:2} that
  \begin{equation}\label{eq:3}
   [a_k,x_1,\dots,x_s]=[a_k,x_{i_1},\dots,x_{i_s}]
  \end{equation}
  if $1\leq k\leq m$, $x_1,\dots,x_s \in \{t_1,\dots, t_n,u_1,\dots,u_n\}$ and $\{i_1,\dots,i_s\}$ is any permutation of 
  $1,\dots,s$. In view of \eqref{eq:1} and \eqref{eq:3}, we see that
  \begin{equation}\label{eq:4}
   [a_k,t_{j_1},t_{j_2},\dots,t_{j_l},u_i]=[a_k,t_{j_1},t_{j_2},\dots,t_{j_l},t_i,t_i]
\end{equation}
if $1\leq k\leq d$, $\{j_1,j_2,\dots,j_l,i\}\subset \{1,2,\dots n\}$. By Lemma \ref{l1} and \eqref{eq:4}, we conclude that
all the elements of $W^+$ can be presented as  linear combinations of the monomials of the following set
$$
 S=\{a_1,\dots, a_m,t_1,\dots,t_m,u_1,\dots,u_m\}\cup \{[a_i,t_{j_1},\dots,t_{j_s}]\mid  i
 ,j_1,\dots,j_s \in\{1,\dots,n\}\}
$$
The product of any two elements of $S$ is defined in the given presentation.
\end{proof}
We will use the following lemmas to show that $W^+$ has a finite presentation.
\begin{lem}\cite[Lemma 5]{bau77}\label{l5}
Let $L$ be a Lie algebra of characteristic $p\neq 2$. Suppose $a,b,t,u$ are elements of $L$ and suppose
$$[a,b]=[a,t,b]=[b,t,a]=[t,u]=0$$
and
$$[a,u]=[a,t,t],\;[b,u]=[b,t,t].$$
Then
$$[[a,\underbrace{t,\dots,t}_{i}],[b,\underbrace{t,\dots,t}_{j}]]=0$$
for every $i\geq 0, j\geq 0$.
 
\end{lem}

\begin{proof}
 Let us denote $a_0=a$, $b_0=b$ and 
 $$a_i=[a,\underbrace{t,\dots,t}_{i}],\;b_j=[b,\underbrace{t,\dots,t}_{j}],\;\text{for}\; i,j\geq 1.$$
 To prove that $[a_i,b_j]=0$ whenever $i,j\geq 0$, we apply induction on $i$ and $j$. For 
 $0\leq i,j\leq 1$, 
 $$[a_0,b_0]=[a_1,b_0]=[a_0,b_1]=0$$
 are given relations. We only need to verify that $[a_1,b_1]=0$ to complete the base cases of induction.
 by Lemma \ref{l0}, $[a_1,b_0]=0$ implies $[a_1,b_1]=[b_0,a_2]$. Similarly, $[b_1,a_0]$ implies
 $[b_1,a_1]=[a_0,b_2]$. So we get
 \begin{equation}\label{l5e1}
  [a_0,b_2]=[a_2,b_0]=[b_1,a_1]
  \end{equation}
In view of Lemma \ref{l0} and the given relations $[a_0,u]=a_2$ and $[b_0,u]=b_2$,
\begin{equation}\label{l5e2}
 [a_0,b_2]=[a_0,[b_0,u]]=[b_0,a_2]
\end{equation}
Combining (\ref{l5e1}) and (\ref{l5e2}), we get $[a_0,b_2]=[b_0,a_2]=[a_2,b_0].$ Since $char(k)\neq 2$,
we conclude that $[a_0,b_2]=0$. Thus $[a_1,b_1]=0$. Now, suppose that
$$[a_i,b_j]=0\;\text{for}\;0\leq i\leq n,\;0\leq j \leq n.$$
Since $[t,u]=0$, by Lemma \ref{l0}, we have 
$$[a_i,u]=a_{i+2},\;[b_i,u]=b_{i+2}\;\text{for any}\; i\in\{1,2,\dots\}.$$
Combining the induction hypothesis with Lemma \ref{l0} , we get
\begin{equation}\label{l5e3}
 [a_i,b_{n+1}]=[b_n,a_{i+1}]=0\;\;\text{for}\; 0\leq i\leq n-1
\end{equation}
and similarly,
\begin{equation}\label{l5e4}
 [b_j,a_{n+1}]=[a_n,b_{j+1}]=0\;\;\text{for}\; 0\leq j\leq n-1
\end{equation}
it remains only to verify that
$$[a_n,b_{n+1}]=[a_{n+1},b_{n+1}]=[b_n,a_{n+1}]= 0.$$
Now $[b_{n-1},a_{n+1}]=0$, so by Lemma \ref{l0}
\begin{equation}\label{l5e5}
[b_{n-1},a_{n+2}]=[b_{n-1},[a_{n+1},t]]=[a_{n+1},b_{n}]
\end{equation}
and $[b_{n-1},a_n]=0$ and $[a_n,b_n]=0$ imply the following relations, respectively.
\begin{equation}\label{l5e6}
[b_{n-1},a_{n+2}]=[b_{n-1},[a_n,u]]=[a_n,b_{n+1}]
 \end{equation}
 \begin{equation}\label{l5e7}
  [a_n,b_{n+1}]=[b_n,a_{n+1}]=-[a_{n+1},b_n]
 \end{equation}
 Putting (\ref{l5e5}), (\ref{l5e6}) and (\ref{l5e7}) together, we get $-[a_{n+1},b_n]=[a_{n+1},b_n]$.
 Since $char(k)\neq 2$, this implies
 \begin{equation}\label{l5e8}
  [a_{n+1},b_n]=0
 \end{equation}
 and a similar argument shows that
 \begin{equation}\label{l5e9}
  [a_{n},b_{n+1}]=0
 \end{equation}
We also need to verify that $[a_{n+1},b_{n+1}]=0$. By Lemma \ref{l0} ,
\begin{equation}\label{l5e10}
 [a_n,b_{n+2}]=[b_{n+1},a_{n+1}]
\end{equation}

\begin{equation}\label{l5e11}
 [a_n,b_{n+2}]=[b_{n},a_{n+2}]
\end{equation}

\begin{equation}\label{l5e12}
 [a_{n+1},b_{n+1}]=[b_{n},a_{n+2}]
\end{equation}
Combining (\ref{l5e10}), (\ref{l5e11}) and (\ref{l5e12}), we get 
$[a_{n+1},b_{n+1}]=[b_{n+2},a_n]=[a_{n+2},b_{n}]=-[a_{n+1},b_{n_1}]$
Therefore,
\begin{equation}\label{l5e13}
 [a_{n+1},b_{n+1}]=0
\end{equation}
Equations (\ref{l5e3}), (\ref{l5e4}), (\ref{l5e8}), (\ref{l5e9}) and (\ref{l5e13}) completes the induction and the proof of 
Lemma \ref{l5}.

\end{proof}

\begin{lem}\label{l6} For any $k,l\in \{1,2,\dots,m\}$
and $ i_1,i_2,\dots,i_s \in \{1,2,\dots,n\}$,
$$[a_k,t_{i_1},t_{i_2},\dots,t_{i_s},a_l]=0$$ 
implies that for any $r\in \{1,\dots,s-1\}$,
$$[[a_l,t_{i_1},t_{i_2},\dots,t_{i_r}],
 [a_k,t_{i_{r+1}},\dots,t_{i_s}]]=0 .$$
 \end{lem}
\begin{proof}By Lemma \ref{l0} and equation \eqref{eq:3} we have
 
  $$\di \begin{array}{lll}
 
[a_k,t_{i_1},t_{i_2},\dots,t_{i_s},a_l]&=&[a_k,t_{i_2},t_{i_3},\dots,t_{i_s},t_{i_1},a_l]\\
 &=&[[a_l,t_{i_1}][a_k,t_{i_2},t_{i_3},\dots,t_{i_s}]\\
 &=&[[a_l,t_{i_1}][a_k,t_{i_3},\dots,t_{i_s},t_{i_2}]\\
 &=&[[a_k,t_{i_3},\dots,t_{i_s}],[a_l,t_{i_1},t_{i_2}]\\
 &\dots&\\
 &=&[[a_l,t_{i_1},t_{i_2},\dots,t_{i_r}],[a_k,t_{i_{r+1}},\dots,t_{i_s}]]\\
 &=& 0.
 \end{array}$$
 \end{proof}

\begin{pro}\label{p2}
  $W^+$ can be presented by the generators
 $$a_1,\dots,a_m,t_1,\dots,t_n,u_1,\dots,u_n,$$
 subject to the finitely many relations
 $$[a_k,t_{j_1},t_{j_2},\dots,t_{j_s},a_l]=0\;\;\;(k,l\in\{1,2,\dots,m\}, 1\leq j_1 <
 j_2<\dots j_s\leq n),$$
 $$[t_i,t_j]=[t_i,u_j]=[u_i,u_j]=0\;\;\;(1\leq i\leq n,\; 1\leq j\leq n),$$
 $$[a_k,u_l]=[a_k,t_l,t_l]\;\;\;(1\leq k\leq n,\; 1\leq l\leq n).$$
\end{pro}
\begin{proof}
 By Lemma \ref{l2}, we only need to show that 
 \begin{equation}\label{eq:5}
  [[a_k,t_{i_1},\dots,t_{i_r}],[a_l,t_{j_1},\dots,t_{j_s}]]=0,
  \end{equation}
  where $1\leq k\leq m,\;1\leq l\leq m,\;i_1,\dots,i_r,j_1,\dots,j_s\in \{1,2,\dots,n\},\;r\geq 0
  ,\;s\geq 0$. To prove this, we apply induction on $r>0$, $s>0$:
  \\If $r=s=1$ we have
  $$w=[[a_k,t_{i_1}],[a_l,t_{j_1}]]$$
  If $i_1=j_1$, we can apply Lemma 5 by taking $a=a_k$, $b=a_l$, $t=t_{i_1}=t_{j_1}$, $u=u_{i_1}$.
  Since we have the relations
  $$[a,b]=[a,t,b]=[b,t,a]=[t,u]=0$$
  $$[a,u]=[a,t,t]\;\text{and}\;[b,u]=[b,t,t],$$
  we get
  $$[[a,t],[b,t]]=0.$$
  If $i_1\neq j_1$, the equation $w=0$ follows from Lemma \ref{l6}.
  Now assume that for $r\leq q$, $s\leq q$, we have
   $$[[a_k,t_{i_1},\dots,t_{i_r}],[a_l,t_{j_1},\dots,t_{j_s}]]=0$$
   Consider
 $$w=[[a_k,t_{i_1},\dots,t_{i_q},t_{i_{q+1}}],[a_l,t_{j_1},\dots,t_{j_s}]]$$
 for $s\leq q$. If $t_{i_1},\dots,t_{i_q},t_{i_{q+1}},t_{j_1},\dots,t_{j_s}$ are all distinct, $w=0$
 is followed from one of the given relations and Lemma \ref{l6}. If not, we prove $w=0$ as follows:
 By Lemma \ref{l0}
   and equation \eqref{l0}, we can assume without loss of generality that
   $$t_{i_{q+1}}=t_{i_q}=t.$$
   Taking $a=[a_k,t_{i_1},\dots,t_{i_{q-1}}]$, $b=[a_l,t_{j_1},\dots,t_{j_s}]$,
   $t=t_{i_q}$ and  $u=u_{i_q}$,  we can apply Lemma 5. Since $a,b,t,u$ satisfy the relations in Lemma \ref{l5},
   we have
   \begin{equation}\label{eq:6}
   w=[[a_k,t_{i_1},\dots,t_{i_q},t_{i_{q+1}}],[a_l,t_{j_1},\dots,t_{j_s}]]=[[a,t,t],b]=0
   \end{equation}
   Similarly, one can show that
   \begin{equation}\label{eq:7}
   [[a_k,t_{i_1},\dots,t_{i_r}],[a_l,t_{j_1},\dots,t_{j_q},t_{j_{q+1}}]]=0
   \end{equation}
   for $r\leq q$.
   To complete the proof, we only need to show that
   
   $$[[a_k,t_{i_1},\dots,t_{i_q},t_{i_{q+1}}],[a_l,t_{j_1},\dots,t_{j_q},t_{j_{q+1}}]]=0$$
   
   and it follows from the equations \eqref{eq:6}, \eqref{eq:7} and Lemma \ref{l5}.
  
 \end{proof}

\section{Proof of Theorem 1} 
Let $W$ be the wreath product of abelian Lie algebras $A$ and $T$ generated by $\{a_1,\dots,a_d\}$ and
$\{t_1,\dots,t_d\}$, respectively. In the previous section we have shown that $W$ is a subalgebra of
a finitely presented Lie algebra $W^+$ generated by $\{a_1,\dots, a_d,t_1,\dots,t_d,u_1,\dots,u_d\}$.
In order to prove
Theorem 1, we compute the growth rate of $W^+$ in this section. In Corollary \ref{c2}, we have shown that
the free metabelian Lie algebra $M$ generated by $d$ elements can be embedded in $W$
 So we
have 
$$\gamma_M \sim n^{d}\precsim \gamma_{W^+}$$
To find an upper bound for the growth rate $\gamma_{W^+}$ of $W^+$, we consider the number of the non-zero monomials 
in $W^+$. In the proof of Lemma \ref{l2}, we have shown that 
all the elements of $W^+$ can be presented as  linear combinations of the monomials of the following set
$$
 S=\{a_1,\dots, a_m,t_1,\dots,t_m,u_1,\dots,u_m\}\cup \{[a_i,t_{j_2},\dots,t_{j_s}]\;\mid\; 1\leq i\leq n
 ,\;j_1,\dots,j_s \in\{1,\dots,n\}\}
$$ and combining this with equation \eqref{eq:3} we see that, as a vector space $W^+$ has a basis which
is a subset of the
following set:
$$\tilde{S}=\{a_1,\dots, a_d,t_1,\dots,t_d,u_1,\dots,u_d\}\cup \{[a_i,t_{j_i},t_{j_2},\dots,t_{j_s}]\mid
1\leq i\leq d,\;1\leq j_1\leq j_2\leq \dots \leq j_s\leq d\}$$
So the growth function $\gamma_{W^+}(n)$ of $W^+$ is less than or equal to the number of elements of length not
greater than
$n$ in $\tilde{S}$,
$$ \di     \gamma_{W^+}(n)\leq 2d+d+\sum_{s=1}^{n-1}d\binom{s+d-1}{d-1}\sim \sum_{s=1}^{n-1}ds^{d-1}\sim n^d$$
and we conclude that 
$$\gamma_{W^+}(n)\sim n^{d}.$$
To complete the proof, we consider the relation between the growth 
functions of a Lie algebra $L$ and its universal enveloping algebra $U(L)$:
Let $L$ denote a Lie algebra generated over a field $k$
by the finite set $X$ whose elements are linearly independent over $k$ and $X^n$ denote the subspace of $L$
spanned by all monomials of length less than or equal to $n$, as we defined in the first section. We may assume $L\subset 
U(L)$, so that $X$ generates $U$ as an associative algebra. Let $u_1,u_2,\dots$ be an ordered basis of $L$ such that 
$X=\{u_1,\dots,u_{\gamma_L(1)}\}$ and $u_{\gamma_L(n-1)}\dots \gamma_{L(n)} $ is a basis for $X^n/X^{n-1}$, $n\geq 2$.
By Poincar\'e-Birkhoff-Witt Theorem \cite {jacob62}, monomials of the form
 $$u_{i_1}\dots u_{i_r}\;\text{with}\; i_1\leq i_2 \leq \dots \leq i_r$$ form a basis for $U$ and we get the following
 relation:
 \begin{equation}\label{eq:9}
\displaystyle \sum_{n=0}^{\infty} b_n t^n =\prod_{n=1}^{\infty}(1-t^n)^{-a_n}
\end{equation}
where $a_n:=dim(X^n/X^{n-1})$ and $b_n$:=number of monomials of length $n$ in $U(L)$
(\cite{smith76}). The relation between the growth rates of $a_n$ and $b_n$ is given in the following proposition,
the proof of which
 can be found in
 various papers (\cite{ber83}, \cite{pet93}, \cite{gribar00}).
\begin{pro}
\label{pro2}
 If $a_n$ and $b_n$ are related by \eqref{eq:9} and
 $a_n\sim n^d$, then 
 $b_n \sim e^{n^{\frac{d+1}{d+2}}}.$
  
 \end{pro}
So the growth of the universal enveloping algebra $U(W^+)$ of $W^+$ is
$$\gamma_{U(W^+)}(n)\sim e^{n^\frac{d}{d+1}}.$$ 

Since $W^+$  is  finitely presented, so is $U(W^+)$  \cite{ufna}. Hence we can conclude that 
for any integer $d>0$, there exists a finitely presented associative algebra of growth $e^{n^\frac{d}{d+1}}$.

\end{section}

\bibliographystyle{alpha}

\bibliography{metabelian}

\end{document}